\documentclass[11pt]{article}
\usepackage{graphicx,color,epsf,bm}
\usepackage{amsmath,amsfonts,amssymb,amscd,bm}
\usepackage{flushend}
\usepackage{hyperref}
\hypersetup{
    colorlinks=true,
    linkcolor=blue,
    filecolor=magenta,      
    urlcolor=cyan,
}

\renewcommand{\vec}[1]{\mbox{\boldmath$#1$}} 

\newcommand{\dif}{\,\mathrm{d}} 

\usepackage{cite}
\title{Semi-Analytical Solution for a Multi-Objective TEAM Benchmark Problem}
\author{Pavel Karban, David P\'{a}nek, Tam\'{a}s Orosz, Ivo Dole\v{z}el \\ University of West Bohemia \\ \{karban, panek50, tamas, idolezel\}@kte.zcu.cz}

\begin{document}
\maketitle

\begin{abstract}
Benchmarking is essential for testing new numerical analysis codes. Their solution is crucial both for testing the partial differential equation solvers and both for the optimization methods. Especially, nature-inspired optimization algorithm-based solvers, where is an important study is to use benchmark functions to test how the new algorithm may perform, in comparison with other algorithms or fine-tune the optimizer parameters. This paper proposes a novel semi-analytical solution of the multi-objective T.E.A.M benchmark problem. The goal of the benchmark problem is to optimize the layout of a coil and provide a uniform magnetic field in the given region.  The proposed methodology was realized in the open-source robust design optimization framework Ārtap, and the precision of the solution is compared with the result of a fully hp-adaptive numerical solver: Agros-suite. The coil layout optimization was performed by derivative-free non-linear methods and the NSGA-II algorithm.
\end{abstract}

\let\thefootnote\relax\footnotetext{Accepted in Periodica Polytechnica Electrical Engineering and Computer Science, the paper will be published in https://pp.bme.hu/eecs}
\section{Introduction}

This paper deals with a seemingly simple problem, where the radius of a given number of circular turns should be optimized that would generate a uniform magnetic field in the prescribed region (Fig.\,\ref{fig:geometry}). This task forms a proposal for the multi-objective Testing Electromagnetic Analysis Methods (TEAM) benchmark problem for Pareto-optimal electromagnetic devices \cite{1,2,3}. TEAM problems offer a wide variety of test problems to benchmark the new partial differential equation and numerical solvers \cite{1,2,3,4,5,6,7}. Moreover, these problems are openly accessible from the website of COMPUMAG society \cite{8}. 

This test problem is inspired by a bio-electromagnetic application for Magnetic Fluid Hyperthermia (MFH), where the uniform magnetic field is used to compare the magnetic properties of the different nanofluids \cite{9,10,11}.
The solenoid design has great importance, as a wide range of application, fields exist, starting from electric power applications \cite{12,13,14,15,16,17,18,19,20,21,22,23,24,25,26,27,28,29,30,31,32,33,34}, through induction heating processes \cite{12}, \cite{24}, \cite{34,35,36,37,38,39,40}⁠ to other biomedical applications \cite{1,2,3} in the industry. The motivation behind the development of the Artap framework \cite{41,42,43}⁠ was a similar problem. To facilitate the development of the induction brazing process, where the main design question was to find a robust design optimum, where the sensitivity of the inductor shape to the  manufacturing tolerances and realized controller design had to be considered together.

The following chapters of the paper propose a novel, semi-analytical solution for the dc uniform magnetic field design and description of the automatized FEM solution using the application of the Ārtap framework.The solution of this benchmark problem is twofold: validating the correctness of the results and demonstrating the applicability of the Ārtap framework.

\section{Problem Description}

\begin{figure}[!ht]
\centering
\includegraphics[width=0.65\linewidth]{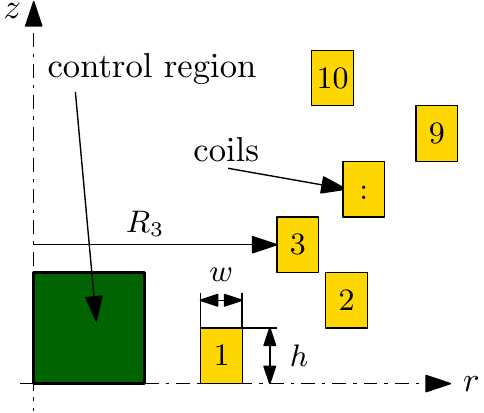} 
\caption{2D axisymmetric model from the upper half the coil geometry and the design variables.}
\label{fig:geometry}
\end{figure}

The described problem is shown in Fig.\,\ref{fig:geometry}.  The task is to get a region (green color in Fig.,\ref{fig:geometry}) with a highly uniform magnetic field distribution. This magnetic field is generated by a prescribed number ($N=20$) of massive circular turns of rectangular cross-section (yellow color in Fig.\,\ref{fig:geometry}).  During the solution of the problem, we are considering only the symmetrical solutions of the problem, as other authors handled this problem  \cite{1,2,3}.

While the dimensions of the turns and the variation of their positions in the z-direction are fixed. The height and the width parameters of the modeled conductors are 1.5 mm and 1.0 mm during the calculations. The inner radius of the turns (radii) can be varied from 5 mm to 50 mm in the $r$-direction. All turns carry a direct current of value $I$, the current density in the conductors are $J_\phi = 2 A/mm^2$. The width and the height of the controlled region are 5 mm in the $r$ and the $z$ directions. In this paper, the following two-goal function based multi-objective version of the task is resolved:
 
\begin{equation}
    F_1(r) = sup_{q=1,np}|\vec{B}(r_q,z_q) - \vec{B_0}(r_q,z_q)|,
\end{equation}

\begin{equation}
    F_2(r) = \Sigma_n R(r_q),
\end{equation}
 
\noindent  where $\vec{B_0} = 2 mT$ is the aimed magnetic flux density, and $\vec{B}$ is the distribution of the magnetic flux density in the field of interest. The value of $\vec{B}$ is calculated in np different points of the region of interest. $F_2$ function represents the mass of the winding, where $n$ is the number of turns, and $R$ is a mass function, which depends on the radii.
 
\section{Semi-analytical Solution}

\begin{figure}[!ht]
\centering
\includegraphics[width=0.92\linewidth]{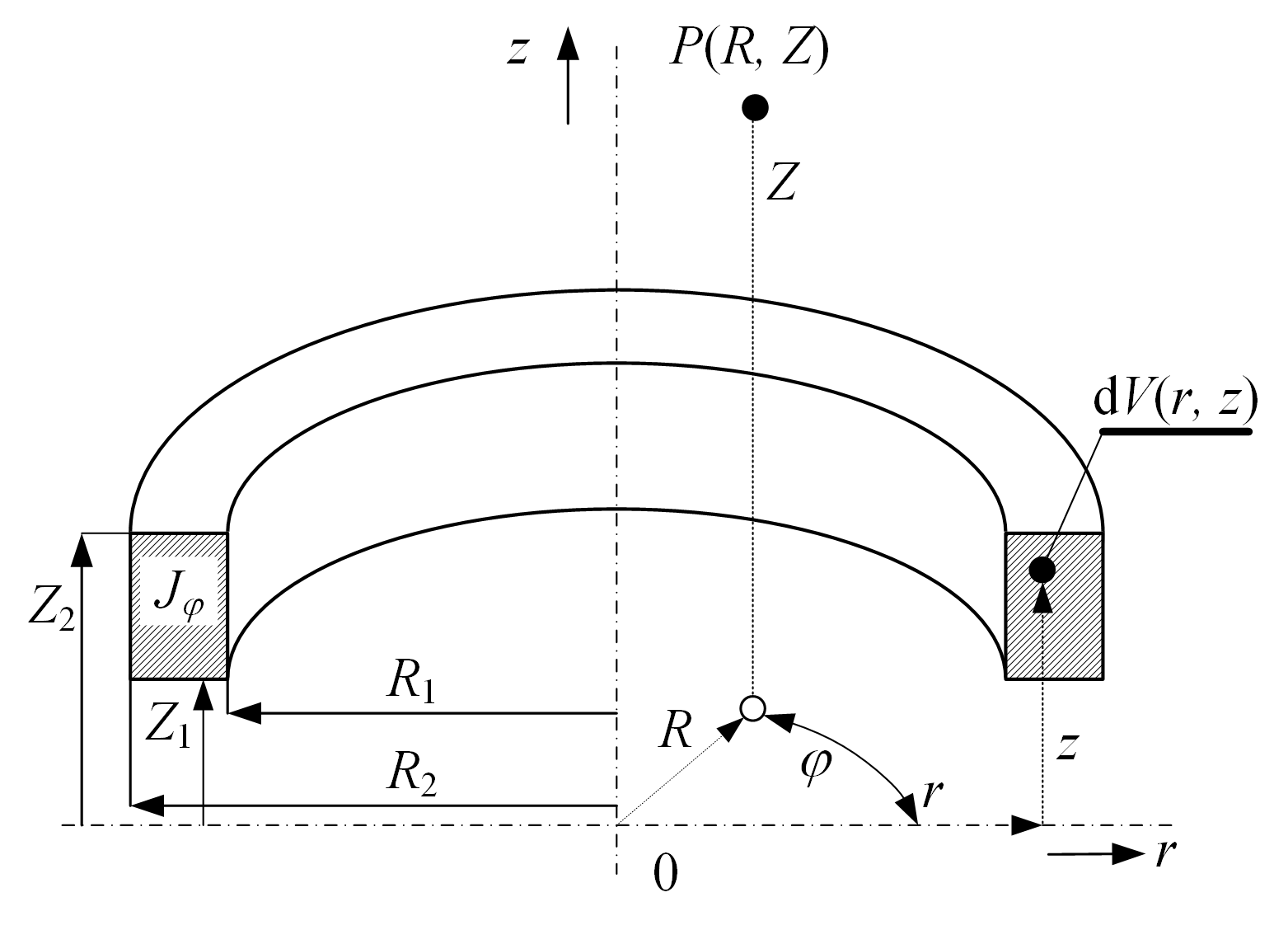} 
\caption{One, massive circular turn of rectangular cross section from the modeled solenoid coil.}
\label{fig:turn}
\end{figure}

Consider a single massive circular turn of rectangular cross section, see Fig.\,\ref{fig:turn}. The turn carries direct current $I$. 

Its density having only the circumferential component $J_{\varphi}$ inversely proportional to the corresponding radius $r$ is given by the formula 
\begin{equation}
J_{\varphi}=\frac{I}{r\,(Z_2-Z_1)\ln(R_2/R_1)}\,.   
\end{equation}

The basic quantity to start with is the circumferential component $A_{\varphi}(R,Z)$ of the magnetic vector potential $\vec A$ that is given by the formula
\begin{equation}
A_{\varphi}(R,Z)=\frac{\mu_0}{4\pi}\int_V\frac{J_{\varphi}
\cos\varphi}{l}\,\dif V\,,   
\end{equation}
where (see Fig.\,2) $l=\sqrt{r^2+R^2-2rR\cos\varphi+(Z-z)^2}$ and $\dif V=r\,\dif r\,\dif\varphi\,\dif z$, $V$ denoting the volume of the ring.

The components of the magnetic flux density $B_r(R,Z)$ and $B_z(R,Z)$ at point $P(R,\,Z)$ in the directions $r$ and $z$ then follow from the relations
\begin{equation}
B_{r}(R,Z)=-\frac{\partial A_{\varphi}(R,Z)}{\partial Z}=
\frac{\mu_0}{4\pi}\int_V\frac{J_{\varphi}(Z-z)\cos\varphi}
{l^3}\,\dif V   
\end{equation}
and
$$
B_{z}(R,Z)=\frac{1}{R}\frac{\partial[R\,A_{\varphi}(R,Z)]}{\partial R}=
$$
\begin{equation}
=\frac{\mu_0}{4\pi R}\int_V\frac{J_{\varphi}[r(r-R\cos\varphi)+
(Z-z)^2]\cos\varphi}{l^3}\,\dif V\,.   
\end{equation}

After substituting for ${J_{\varphi}}$ from (1), for $l$ and for $\dif V$, we obtain
\begin{equation}
B_{r}(R,Z)=C\int_{\varphi=0}^{2\pi}\int_{r=R_1}^{R_2}\int_{z=Z_1}^{Z_2}I_1\,\dif z\,\dif r\,\dif \varphi\,,   
\end{equation}
where
$$
I_1=\frac{(Z-z)\cos\varphi}{\sqrt{[r^2+R^2-2rR\cos\varphi+(Z-z)^2]^3}}
$$
and 
\begin{equation}
B_{z}(R,Z)=\frac{C}{R}\int_{\varphi=0}^{2\pi}\int_{r=R_1}^{R_2}\int_{z=Z_1}^{Z_2}I_2\,\dif z\,\dif r\,\dif \varphi\,,   
\end{equation}
where
$$
I_2=\frac{[r(r-R\cos\varphi)+(Z-z)^2]\cos\varphi}{\sqrt{[r^2+R^2-2rR\cos\varphi+(Z-z)^2]^3}}\,.
$$
The constant $C$ occurring in both (5) and (6) is given as
\begin{equation}
C=\frac{\mu_0 I}{4\pi (Z_2-Z_1)\ln(R_2/R_1)}\,.
\end{equation}

After the double integration with respect to $r$ and $z$, the components of the magnetic flux density are given by the formulae \cite{dolezel}
$$
B_{r}(R,Z)=C\cdot[g(R_2,R,Z_2-Z)-g(R_2,R,Z_1-Z)-
$$
\begin{equation}
-g(R_1,R,Z_2-Z)+g(R_1,R,Z_1-Z)]
\end{equation}
and
$$
B_{z}(R,Z)=C\cdot[h(R_2,R,Z_2-Z)-h(R_2,R,Z_1-Z)-
$$
\begin{equation}
-h(R_1,R,Z_2-Z)+h(R_1,R,Z_1-Z)]\,.
\end{equation}

Here, for example
\begin{equation}
g(R_2,R,Z_2-Z)=\int_{\varphi=0}^{2\pi}\ln[R_2-R\cdot\cos\varphi+d_{22}]\cos\varphi\,{\rm d}\varphi\,,
\end{equation}
\begin{equation}
h(R_2,R,Z_2-Z)=-\int_{\varphi=0}^{2\pi}\ln[Z_2-Z+d_{22}]\,{\rm d}\varphi\,,
\end{equation}
and
\begin{equation}
d_{22}=\sqrt{R_2^2+R^2-2R_2 R\cos\varphi+(Z-Z_2)^2}\,.
\end{equation}

The other functions are obtained by standard interchanging of the indices. The last integrals with respect to $\phi$ are calculated using the Gauss quadrature formulae. Magnetic field produced by more turns is then given by the superposition of the partial fields produced by particular turns. The computations were realized in the Ārtap framework. It is the part of the package, which can be downloaded from the homepage of the project:
\url{http://pypi.org/project/artap} .

The accuracy of the above results is compared with the results of Agros \cite{43}.

\section{FEM Model}

\begin{figure}[!ht]
\centering
\includegraphics[width=0.85\linewidth]{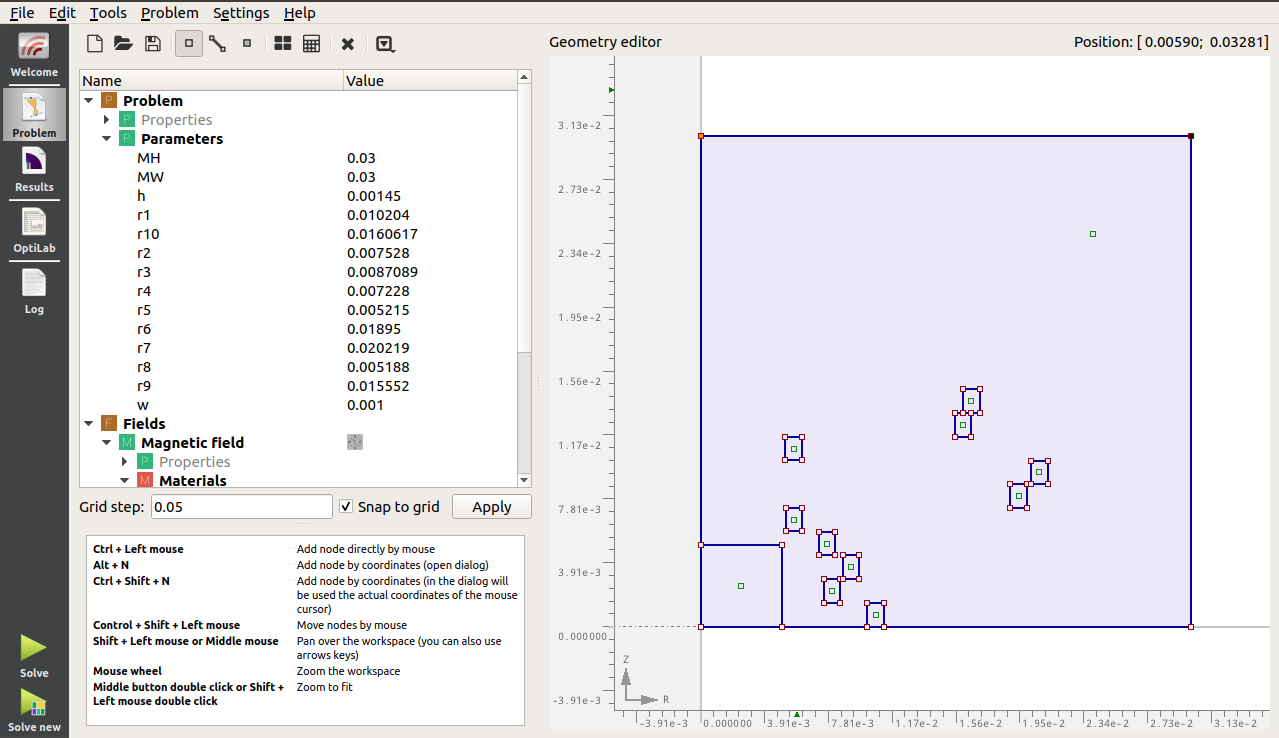} 
\caption{The realized parametric geometry for the TEAM benchmark problem in Agros.}
\label{fig:agros_ui}
\end{figure}

A precompiled version of the fully hp-adaptive FEM solver: Agros Suite is integrated into the Ārtap framework, and it can be invoked via a Python scripting language \cite{41,42}. The model geometry can be made by the aid of the user Fig. \ref{fig:agros_ui} The realized parametric geometry for the TEAM benchmark problem in Agros interface. 

Then we can use the ‘Create script from the model’ function of the Problem toolbar to copy the tested FEM model into the Python script (Fig. \ref{fig:agros_ui}). This script can be inserted directly into the optimization code. The modeler has to connect the Agros model parameters with the model parameters of the Ārtap project at the beginning of the code and its ready to use.

\begin{figure}[!ht]
\centering
\includegraphics[width=0.75\linewidth]{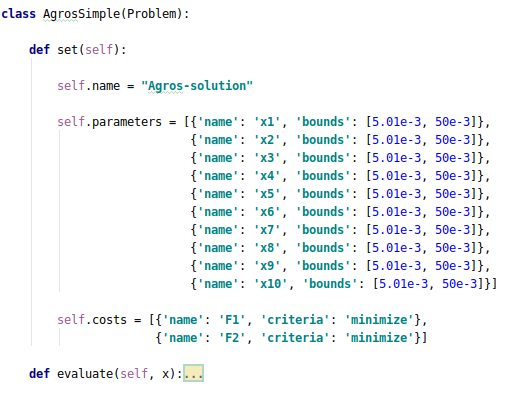} 
\caption{Initialization of the Optimization problem in Ārtap.}
\label{fig:code}
\end{figure}

The optimization task can be initialized by the usage of the \textit{set()} function (Fig \ref{fig:code}). Here, four parameters and methods have to be defined: 

\begin{itemize}
    \item the name of the optimization task.
    \item a dictionary list with the optimized parameters.
    \item the objective functions.
\end{itemize}

The exported python code has to be inserted into the \textit{evaluate()} method. It gets the optimization variables via the parameter vector (\textit{x.vector}) of the Individual class. This parameter introduced to the code to give the possibility to add other non-optimized, FEM or analytical results of the calculation to the Individual. These values can be saved into a database and can be post-processed after the calculation.
The database management and the parallelization can be made automatically; it only has to be defined by a keyword. In this paper, we are using the NSGA-II \cite{41,44} algorithm to solve the optimization task. This optimization method does not need initial values, so only the name of the parameter and the lower and the upper bounds are enough for the initialization. The NSGAII algorithm is contained by the algorithm class. It is initialized by a maximum of 100 individuals in 100 generations. The realized problem can be downloaded from the project homepage, and it is part of the python package \url{https://pypi.org/project/artap/}.

\section{Results and Discussion}

\begin{figure}[!htb]
\centering
\includegraphics[width=0.85\linewidth]{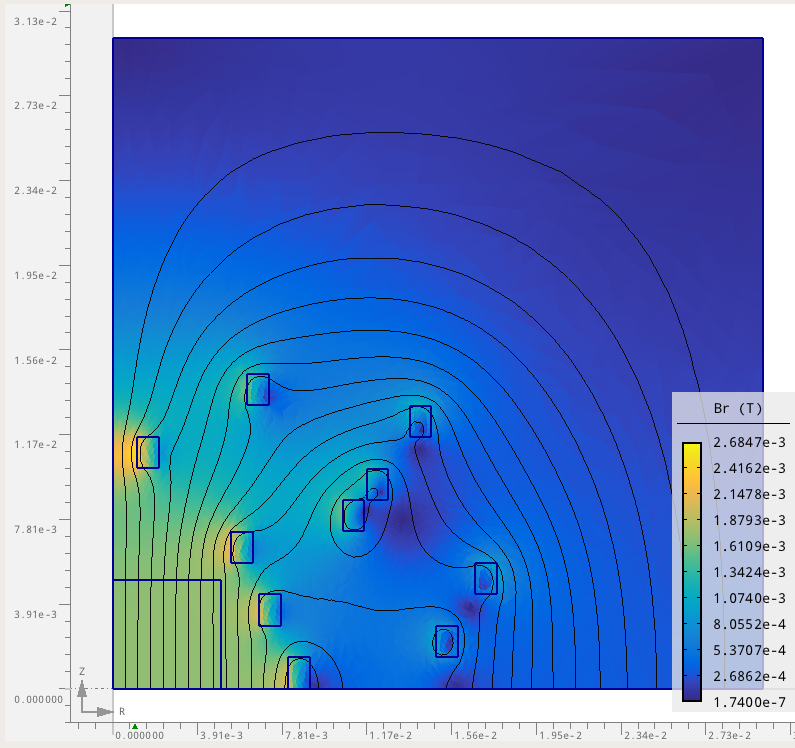} 
\caption{The examined geometry for the TEAM benchmark problem and the flux density distribution in Agros.}
\label{fig:femresults}
\end{figure}

The results of the semi-analytical calculation compared with a FEM calculation to benchmark and validate its
results. For comparison, one possible turn layout is selected. Here, the radii of the turns set by the list of the x parameters: $x= [0.00808, 0.0149, 0.00674, 0.0167, 0.00545, 0.0106, 0.0117, 0.0111, \newline 0.01369, 0.00619]$. 
Where every value is given in m, this solenoid layout and the resulting flux density distribution is depicted in Fig \ref{fig:femresults}. 

\begin{figure}[!htb]
\centering
\includegraphics[width=0.85\linewidth]{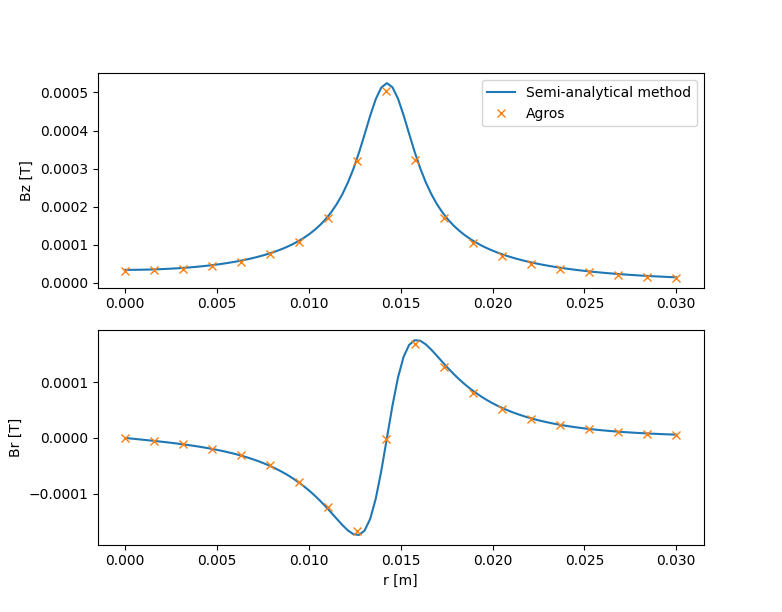} 
\caption{The $B_r$ and $B_z$ flux density components along the $x$ = 0.003 m vertical line, the values were calculated in 20 points.}
\label{fig:analytic}
\end{figure}

The radial ($B_r$) and the axial ($B_z$) components of the magnetic flux density were compared along a vertical line ($r = 0.003$), where 10 different points are selected from the area of interest. The results are compared in Fig \ref{fig:analytic}.

\begin{figure}[!htb]
\centering
\includegraphics[width=0.85\linewidth]{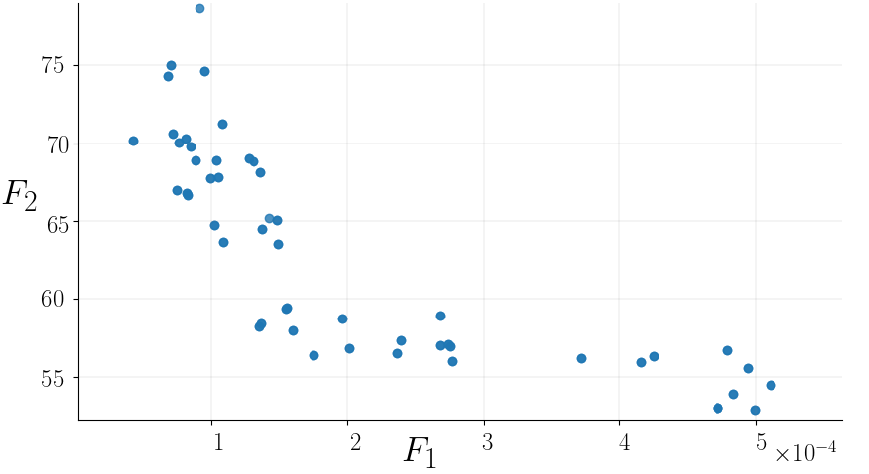} 
\caption{Last generation of the optimization for $F_1$ and $F_2$.}
\label{fig:optres}
\end{figure}

As it can be seen from the results, the difference between the FEM and the semi-analytical solution is negligible. However, the computational cost of the semi-analytical calculation is much lesser than a single FEM solution. Therefore, to accelerate the optimization process, this formulation is used to search the Pareto-solutions (Fig. \ref{fig:optres}). During the calculation, 10000 iterations were performed. The shape of the resulting Pareto-front is similar to the solution, which is presented in the proposal of the benchmark problem \cite{1,2,3}.

\section{Conclusions}

A novel semi-analytical solution was proposed to the multi-objective TEAM benchmark problem for electromagnetic devices. The proposed methodology was realized in the open-source robust design optimization framework Ārtap, and the precision of the solution is compared with the result of a fully hp-adaptive numerical solver: Agros-suite. The coil layout optimization was performed by derivative free non-linear optimization method, the NSGA-II algorithm. The provided semi-analytical solution can be used to test other FEM solvers and significantly reduces the calculation time of the optimization process.

\section*{Acknowledgement}
This research has been supported by the Ministry of Education, Youth and Sports of the Czech Republic under the RICE New Technologies and Concepts for Smart Industrial Systems, project no. LO1607 and by an internal project SGS-2018-043.

\bibliographystyle{IEEEtran}
\bibliography{references}
\end{document}